\documentclass{article}%
\usepackage{amsmath}
\usepackage{amsfonts}
\usepackage{amssymb}
\usepackage{graphicx}%
\providecommand{\U}[1]{\protect \rule{.1in}{.1in}}
\newtheorem{theorem}{Theorem}

\newtheorem{lem}[theorem]{Lemma}

\newtheorem{remark}[theorem]{Remark}

\newenvironment{proof}[1][Proof]{\noindent \textbf{#1.} }{\  \rule{0.5em}{0.5em}}

\newcommand{\abs}[1]{\left\vert#1\right\vert}
\newcommand{\set}[1]{\left\{#1\right\}}
\newcommand{\Real}{\mathbb R}
\newcommand{\eps}{\varepsilon}
\newcommand{\lm}{\lambda}

\newcommand{\F}{\mathcal{F}}

\newcommand{\Hh}{\mathcal{H}^2}
\newcommand{\Ss}{\mathcal{S}^2}

\newcommand{\ul}{\underline}

\begin{document}

\title{The Equivalence between Uniqueness and Continuous Dependence of Solution for BSDEs with Continuous Coefficient}
\author{Guangyan JIA$^a$ \hspace*{.15cm} Zhiyong YU$^{a,b}$\thanks{Corresponding author. Email address: yuzhiyong@sdu.edu.cn}\\
{\small $^a$ School of Mathematics and System Sciences, Shandong
University, Jinan  250100, China}\\
{\small $^b$School of Economics, Shandong University, Jinan 250100,
China.}}
\date{}
\maketitle

\begin{abstract}
In this paper, we will prove that, if the coefficient  $g=g(t,y,z)$
of a BSDE is assumed to be continuous and linear growth in $(y,z)$,
then the uniqueness of solution and continuous dependence with
respect to  $g$ and the terminal value $\xi$ are equivalent.
\end{abstract}



\section{Introduction}

We consider the following $1$--dimensional backward stochastic differential
equation (BSDE):%
\begin{equation}
y_{t}=\xi+\int_{t}^{T}g(s,y_{s},z_{s})\,ds-\int_{t}^{T}z_{s}\,dW_{s},\qquad
t\in \lbrack0,T].\label{bsdeequ}%
\end{equation}
where the terminal condition $\xi$ and the coefficient $g=g(t,y,z)$ are given.
$W$ is a $d$--dimensional Brownian motion. The solution $(y_{t},z_{t}%
)_{t\in \lbrack0,T]}$ is a pair of square integrable processes. A
foundational and interesting problem is: what is the relationship
between the uniqueness of solution and continuous dependence with
respect to $g$ or $\xi$? In the standard situation where $g$
satisfies linear growth condition and Lipschitz condition in
$(y,z)$, it has been proved by Pardoux and Peng \cite{PP1990} that
there exists a unique solution. In this case, the continuous
dependence with respect to $g$ and $\xi$ is is described by the
following inequality (see El Karoui, Peng and Quenez
\cite{EPQ1997}):
\begin{equation}\label{Lip-CD}
E\left\{\sup_{0\leq t\leq T}\abs{y^1_t-y^2_t}^2\right\} \leq C
E\left\{\abs{\xi^1-\xi^2}^2
+\int_0^T\abs{g^1(t,y^1_t,z^1_t)-g^2(t,y^1_t,z^1_t)}^2dt \right\},
\end{equation}
where $(y^1_{t},z^1_{t})_{t\in \lbrack0,T]}$ and
$(y^2_{t},z^2_{t})_{t\in \lbrack0,T]}$ are the unique solutions of
BSDE $(g^1,\xi^1)$ and BSDE $(g^2,\xi^2)$ respectively. From this,
fruitful results are derived. However in the case where $g$ is only
continuous in $(y,z)$, in place of the Lipschitz condition,
Lepeltier and San Martin \cite{LepeltierMartin1997} have proved that
there is at least one solution. In fact, there is either one or
uncountable many solutions in this situation(see Jia and Peng
\cite{JP2006}). To answer the question whether the uniqueness of
solution also implies the continuous dependence with respect to $g$
and $\xi$ is what this paper will achieve.

In this paper we will prove that if the coefficient $g$ satisfies
the conditions given in \cite{LepeltierMartin1997}, then the
uniqueness of solution and continuous dependence with respect to $g$
and $\xi$ are equivalent. This result, which can be regarded as the
analog of the inequality (\ref{Lip-CD}) in some sense, provides a
useful method to study BSDEs with continuous coefficient.

This paper is organized as follows. In Section \ref{secPre} we
formulate the problem accurately and give some preliminary results.
Section \ref{secMainResults} is devoted to proving the equivalence
between uniqueness and continuous dependence with respect to
terminal value $\xi$. Finally, in Section \ref{secGeneralcase} we
will prove the equivalence of uniqueness and continuous dependence
with respect to parameters $g$ and $\xi$.

\section{Preliminaries}\label{secPre}

Let $(\Omega,\mathcal{F},P)$ be a probability space and
$(W_{t})_{t\geq0}$ be a $d$-dimensional standard Brownian motion in
this space. Let $(\mathcal{F}_{t})_{t\geq0}$ be the filtration
generated by this Brownian motion: $\mathcal{F}_{t}=\sigma \left \{
W_{s},s\in \lbrack 0,t]\right \}  \vee \mathcal{N}$,
$\mathbb{F}=(\mathcal{F}_{t})_{t\geq0}$, where $\mathcal{N}$ is the
set of all $P$-null subsets.

Let $T>0$ be a fixed real number. In this paper, we always work in
the space $(\Omega,\mathcal{F}_{T},P)$. For a positive integer $n$
and $z\in \mathbb{R}^{n}$, we denote by $\left \vert z\right \vert
$ the Euclidean norm of $z$. We will denote by $\mathcal{H}_{n}^{2}=\mathcal{H}_{n}^{2}(0,T;\mathbb{R}%
^{n})$, the space of all $\mathbb{P}$--progressively measurable $\mathbb{R}%
^{n}$--valued processes s.t. $E[\int_{0}^{T}\left \vert \psi
_{t}\right \vert ^{2}\,dt]<\infty$, and by $\mathcal{S}^{2}=\mathcal{S}%
^{2}(0,T;\mathbb{R})$ the elements in
$\mathcal{H}_{n}^{2}(0,T;\mathbb{R})$ with continuous paths s.t.
$E[\sup_{t\in \lbrack0,T]}\left \vert \psi_{t}\right \vert
^{2}]<\infty$.

The coefficient $g$ of BSDE is a function $g(\omega,t,y,z):\Omega \times
\lbrack0,T]\times \mathbb{R}\times \mathbb{R}^{d}\rightarrow \mathbb{R}$
satisfying the following assumptions:

(H1): linear growth: there exists a nonnegative constant $A$, such
that $\left \vert g(\omega ,t,y,z)\right \vert \leq A(1+\left
\vert y\right \vert +\left \vert z\right \vert )$,\ $\forall
t,\omega,y,z$

(H2): $(g(t,y,z))_{t\in \lbrack0,T]}\in \mathcal{H}_{1}^{2}$,\  \ for each
$(y,z)\in \mathbb{R\times R}^{d}$

(H3): $g(\omega,t,.,.)$ is continuous for fixed $(t,\omega)$.

Given by Lepeltier and San Martin \cite[Th. 1]{LepeltierMartin1997},
under (H1)---(H3) and for each given $\xi \in
L^{2}(\Omega,\mathcal{F}_{T},P)$, there exists at least one solution
$(y_{t},z_{t})_{t\in \lbrack0,T]}\in \mathcal{S}^{2}\times
\mathcal{H}_{d}^{2}$ of BSDE (\ref{bsdeequ}).
\cite{LepeltierMartin1997} also gives the existence of the maximal
solution $(\bar{y}_{t},\bar{z}_{t})_{t\in \lbrack0,T]}$ and the
minimal solution $(\underline{y}_{t},\underline{z}_{t})_{t\in
\lbrack0,T]}$ of BSDE (\ref{bsdeequ}) in the sense that any solution
$(y_{t},z_{t})_{t\in \lbrack0,T]}\in \mathcal{S}^{2}\times
\mathcal{H}_{d}^{2}$ of BSDE (\ref{bsdeequ}) must satisfy
$\underline{y}_{t}\leq y_{t}\leq \bar{y}_{t}$, a.s., for all $t\in
\lbrack0,T]$.

It is well known that under the standard assumptions where $g$ is
Lipschitz continuous in $(y,z)$,  for any random variable $\xi$ in
$L^2({\cal F}_T)$, the BSDE (\ref{bsdeequ}) has a unique adapted
solution, say $(y_t,z_t)_{t\in[0,T]}$ such that $z\in{\cal H}_d^2$
and $y\in{\cal S}^2$ (see \cite{PP1990}). And we have the following
estimate for solution of BSDEs with Lipschitz continuous generator
$g$ coming from \cite{EPQ1997}.

\begin{lem}\label{lem 2.1}
If $\xi^1, \xi^2 \in L^2({\cal F}_T)$ and $g$ is Lipschitz
continuous in $(y,z)$. Then, for the solutions
$(y_t^1,z_t^1)_{t\in[0,T]}$ and $(y_t^2,z_t^2)_{t\in[0,T]}$ of the
BSDEs,$(g,T,\xi^1)$ and $(g,T,\xi^2)$ respectively, we have
$$E[\sup_{0\leq t\leq T}|y_t^1-y_t^2|^2]\leq CE|\xi^1-\xi^2|^2$$
where $C$ is a positive constant only depending on Lipschitz
constant of $g$.
\end{lem}

Now, we will recall some properties and associated approximation
about BSDEs with $g$ satisfying Assumptions (H1)--(H3)(see
\cite{LepeltierMartin1997} for details).
\begin{lem}\label{lem2.2}
If $g$ satisfies Assumptions (H1)---(H3), and we set
$$\ul{g}_m(t,y,z):=\inf_{(u,v)\in\Real^{1+d}}\set{g(t,u,v)+m(\abs{y-u}+\abs{z-v})},$$
and $$\bar{g}_m(t,y,z):=\sup_{(u,v)\in
\Real^{1+d}}\set{g(t,u,v)-m(\abs{y-u}+\abs{z-v})},$$ then for any
$m\ge A$, we have

(1).  For any $y\in\Real$, $z\in\Real^d$ and $t\in[0,T]$,
$\ul{g}_m(t,y,z)\le A(\abs{y}+\abs{z}+1)$, and
$\bar{g}_m(t,y,z)\le A(\abs{y}+\abs{z}+1)$.

(2). For any $y\in\Real$, $z\in\Real^d$ and $t\in[0,T]$,
$\ul{g}_m(t,y,z)$ is non-decreasing in $m$ and $\bar{g}_m(t,y,z)$
is non-increasing in $m$.

(3). $\ul{g}_m$ and $\bar{g}_m$ are Lipschitz functions, i.e., for
any $y_1,y_2\in\Real$, $z_1,z_2\in\Real^d$ and $t\in[0,T]$,
$\abs{\ul{g}_m(t,y_1,z_1)-\ul{g}_m(t,y_2,z_2)}\le
m(\abs{y_1-y_2}+\abs{z_1-z_2})$ and
$\abs{\bar{g}_m(t,y_1,z_1)-\bar{g}_m(t,y_2,z_2)}\le
m(\abs{y_1-y_2}+\abs{z_1-z_2})$.

(4). If $(y_m,z_m)\to (y,z)$ as $m\to\infty$, then
$\ul{g}_m(t,y_m,z_m)\to g(t,y,z)$ and $\bar{g}_m(t,y_m,z_m)\to
g(t,y,z)$ as $m\to\infty$.
\end{lem}

\begin{lem}\label{lem2.3}
If the processes $(\ul{y}_t^m,\ul{z}_t^m)_{t\in[0,T]}$ and
$(\bar{y}_t^m,\bar{z}_t^m)_{t\in[0,T]}$ are the unique solutions
of the BSDEs $(\ul{g}_m,T,\xi)$ and $(\bar{g}_m,T,\xi)$
respectively, then
$$(\ul{y}_t^m,\ul{z}_t^m)_{t\in[0,T]}\to(\ul{y}_t,\ul{z}_t)_{t\in[0,T]},\ and\ (\bar{y}_t^m,\bar{z}_t^m)_{t\in[0,T]}\to(\bar{y}_t,\bar{z}_t)_{t\in[0,T]},\ (m\to\infty)$$
in $\Ss\times\Hh_d$, where $(\ul{y}_t,\ul{z}_t)_{t\in[0,T]}$ and
$(\bar{y}_t,\bar{z}_t)_{t\in[0,T]}$ are the minimal solution and
maximal solution of BSDE (\ref{bsdeequ}).
\end{lem}

\section{Main Results}\label{secMainResults}
In this section, we will prove the equivalence of uniqueness of
solution and continuous dependence with respect to terminal value
$\xi$.
\begin{theorem}\label{thm3.1}
If Assumptions (H1)---(H3) hold for $g$,  then the following two
statements are equivalent.

(i). Uniqueness: The equation (\ref{bsdeequ}) has a unique
solution.

(ii). Continuous dependence with respect to $\xi$: For any
$\set{\xi_n}_{n=1}^{\infty}$, $\xi\in L^2({\cal F}_T)$, if
$\xi_n\to\xi$ in $L^2({\cal F}_T)$ as $n\to\infty$, then
\begin{equation}\label{3.1}
\lim_{n\to\infty}E[\sup_{t\in[0,T]}\abs{y_t^{\xi_n}-y_t^{\xi}}^2]=0
\end{equation}
where $(y_t^{\xi},z_t^{\xi})_{t\in[0,T]}$ is any solution of BSDE
(\ref{bsdeequ}) and $(y_t^{\xi_n},z_t^{\xi_n})_{t\in[0,T]}$ are
any solutions of the BSDEs $(g,T,\xi^n)$.
\end{theorem}
\begin{proof}
Firstly, we will prove that (i) implies (ii). Given $n$, we note
that for any solution $(y^{\xi_n}_t,z^{\xi_n}_t)_{t\in[0,T]}$ of
BSDE $(g,T,\xi^n)$, we have
\begin{equation}\label{3.3}
\ul{y}_t^{\xi_n}\le y_t^{\xi_n}\le \bar{y}_t^{\xi_n},P-a.s.\qquad
t\in[0,T],
\end{equation}
Now, we consider the following equations:
\begin{equation}\label{3.4}
\ul{y}_t^{m,\xi_n}=\xi_n+\int_t^T
\ul{g}_m(s,\ul{y}_s^{m,\xi_n},\ul{z}_s^{m,\xi_n})\,ds -\int_t^T
\ul{z}_s^{m,\xi_n}\,dW_s
\end{equation}
and
\begin{equation}\label{3.5}
\bar{y}_t^{m,\xi_n}=\xi_n+\int_t^T
\bar{g}_m(s,\bar{y}_s^{m,\xi_n},\bar{z}_s^{m,\xi_n})\,ds -\int_t^T
\bar{z}_s^{m,\xi_n}\,dW_s
\end{equation}
where $(\ul{y}_t^{m,\xi_n},\ul{z}_t^{m,\xi_n})_{t\in[0,T]}$ and
$(\bar{y}_t^{m,\xi_n},\bar{z}_t^{m,\xi_n})_{t\in[0,T]}$ are unique
solutions of (\ref{3.4}) and (\ref{3.5}) respectively.

Thanks to Lemma \ref{lem2.3}, we know that
$$(\ul{y}_t^{m,\xi_n},\ul{z}_t^{m,\xi_n})\to(\ul{y}_t^{\xi_n},\ul{z}_t^{\xi_n}),\ and\ (\bar{y}_t^{m,\xi_n},\bar{z}_t^{m,\xi_n})\to(\bar{y}_t^{\xi_n},\bar{z}_t^{\xi_n}),\ t\in[0,T].$$
in $\Ss\times\Hh_d$ as $m\to\infty$, and get the following
inequalities
\begin{equation}\label{3.6}
\ul{y}_t^{m,\xi_n}\le \ul{y}_t^{\xi_n}\le y_t^{\xi_n}\le
\bar{y}_t^{\xi_n}\le \bar{y}_t^{m,\xi_n},\quad for\ any\ n,
t\in[0,T]\ and\ m\ge A\
\end{equation}

From inequality (\ref{3.6}), we have
\begin{eqnarray*}
y_t^{\xi_n}-y_t^\xi&=&y_t^{\xi_n}-\bar{y}_t^{m,\xi_n}+\bar{y}_t^{m,\xi_n}-\bar{y}_t^{m,\xi}+\bar{y}_t^{m,\xi}-y_t^\xi\\
&\le&(\bar{y}_t^{m,\xi_n}-\bar{y}_t^{m,\xi})+(\bar{y}_t^{m,\xi}-y_t^\xi)
\end{eqnarray*}
and
\begin{eqnarray*}
y_t^{\xi_n}-y_t^\xi&=&y_t^{\xi_n}-\ul{y}_t^{m,\xi_n}+\ul{y}_t^{m,\xi_n}-\ul{y}_t^{m,\xi}+\ul{y}_t^{m,\xi}-y_t^\xi\\
&\ge&(\ul{y}_t^{m,\xi_n}-\ul{y}_t^{m,\xi})+(\ul{y}_t^{m,\xi}-y_t^\xi)
\end{eqnarray*}
Thus
\begin{eqnarray*}
E[\sup_{t\in[0,T]}\abs{y_t^{\xi_n}-y_t^{\xi}}^2] &\le&
2E[\sup_{t\in[0,T]}\abs{\ul{y}_t^{m,\xi_n}-\ul{y}_t^{m,\xi}}^2]
+2E[\sup_{t\in[0,T]}\abs{\ul{y}_t^{m,\xi}-y_t^{\xi}}^2]\\
&+&2E[\sup_{t\in[0,T]}\abs{\bar{y}_t^{m,\xi_n}-\bar{y}_t^{m,\xi}}^2]+2E[\sup_{t\in[0,T]}\abs{\bar{y}_t^{m,\xi}-y_t^{\xi}}^2]\\
\end{eqnarray*}
where $(\ul{y}^{m,\xi}_t,\ul{z}^{m,\xi}_t)_{t\in[0,T]}$ and
$(\bar{y}^{m,\xi}_t,\bar{z}^{m,\xi}_t)_{t\in[0,T]}$ are solutions
of BSDEs $(\ul{g}_m,T,\xi)$ and $(\bar{g}_m,T,\xi)$ respectively.

By Lemma \ref{lem 2.1} and Lemma \ref{lem2.2}, as $n\to\infty$, we
have
$$E[\sup_{t\in[0,T]}\abs{\ul{y}_t^{m,\xi_n}-\ul{y}_t^{m,\xi}}^2]\to 0,\ and\
E[\sup_{t\in[0,T]}\abs{\bar{y}_t^{m,\xi_n}-\bar{y}_t^{m,\xi}}^2]\to
0,\quad for\ any\ m.$$

By Lemma \ref{lem2.3} and the uniqueness of solution for BSDE
(\ref{bsdeequ}), we get
$$
E[\sup_{t\in[0,T]}\abs{\ul{y}_t^{m,\xi}-y_t^{\xi}}^2]\to 0, \ and\
E[\sup_{t\in[0,T]}\abs{\bar{y}_t^{m,\xi}-y_t^{\xi}}^2]\to 0
$$
as $m\to\infty$. That is (ii).

Now, we will prove that (ii) implies (i). We take $\xi_n=\xi$. For
equations $(g,T,\xi^n)$, we set
$y_t^{\xi_n}=\bar{y}_t^{\xi_n}=\bar{y}_t^{\xi}$. For the equation
(\ref{bsdeequ}), we set $y_t^{\xi}=\ul{y}_t^{\xi}$. From (ii), we
have $\bar{y}_t^{\xi}=\ul{y}_t^{\xi}$. The proof is complete.
\end{proof}
\begin{remark}\label{rem3.3}
In fact, when the solution of (\ref{bsdeequ}) is not unique, the
continuous dependence may not hold true in general. For example,
we take $g(t,y,z)=3y^{2/3}$, $\xi=0$. It is easy to know that
$(y_t,z_t)_{t\in[0,T]}=(0,0)_{t\in[0,T]}$ and
$(Y_t,Z_t)_{t\in[0,T]}=((T-t)^3,0)_{t\in[0,T]}$ both are solutions
of BSDE
$$
y_t=\int_t^T 3y^{\frac{2}{3}}_s ds -\int_t^T z_s dW_s;\quad 0\leq
t\leq T.
$$
Set $\xi_n=1/n$, the BSDEs
$$
y_t=\frac{1}{n}+\int_t^T 3y^{\frac{2}{3}}_s ds -\int_t^T z_s
dW_s;\quad 0\leq t\leq T, \quad n=1,2,\cdots.
$$
have unique solutions
$(y^{\frac{1}{n}}_t,z^{\frac{1}{n}}_t)=((T-t+\frac{1}{\sqrt[3]{n}})^3,0)$
for $n=1,2,\cdots$. But
\begin{eqnarray*}
\lim_{n\to\infty}E[\sup_{t\in[0,T]}\abs{y^{\frac{1}{n}}_t-y_t}^2]=T^6\neq
0=\lim_{n\to\infty}E[\sup_{t\in[0,T]}\abs{y^{\frac{1}{n}}_t-Y_t}^2]\\
\end{eqnarray*}
\end{remark}

\section{The General Case}\label{secGeneralcase}
In this section, we will deal with the more general case, that is,
the relationship between uniqueness of solution and continuous
dependence with respect not only to $\xi$ but also to $g$. Now, we
consider the following BSDEs:
\begin{equation}\label{4.1}
y_t^{\lm}=\xi^{\lm}+\int_t^T g^{\lm}(s,y_s^{\lm},z_s^{\lm})\,ds
-\int_t^Tz_s^{\lm}\,dW_s,
\end{equation}
where $\lm$ belongs to a nonempty set $D\subset\Real$. The
coefficient $g^\lm$ is a function $g(\omega,t,y,z):D\times\Omega
\times \lbrack0,T]\times \mathbb{R}\times
\mathbb{R}^{d}\rightarrow \mathbb{R}$ satisfying the following
assumptions:

(H1'): linear growth: there exists a nonnegative constant $A$,
such that $\left \vert g^\lm(\omega ,t,y,z)\right \vert \leq
A(1+\left \vert y\right \vert +\left \vert z\right \vert )$,\
$\forall \lm,t,\omega,y,z$.

(H2'): $(g(t,y,z))_{t\in \lbrack0,T]}\in \mathcal{H}_{1}^{2}$,\  \
for each $(y,z)\in \mathbb{R\times R}^{d}$ and $\lm\in D$.

(H3'): $g(\omega,t,.,.)$ is continuous for fixed $(t,\omega,\lm)$.

(H4'): uniform continuity: $g^{\lm}$ is continuous in $\lm=\lm_0$
uniformly with respect to $(y,z)$.

When (H1') and (H3') are replaced by Lipschitz condition (L),
i.e., there exists a nonnegative constant $K$, such that $\left
\vert g^\lm(\omega ,t,y_1,z_1)-g^\lm(\omega ,t,y_2,z_2)\right
\vert \leq K(\left \vert y_1-y_2\right \vert +\left \vert
z_1-z_2\right \vert )$,\ $\forall \lm,t,\omega,y_1,z_1$ and
$y_2,z_2$, the BSDE (\ref{4.1}) has a unique adapted solution for
any $\lm\in D$. And we have the following property:
\begin{lem}\label{lem 4.1}
If $\xi^{\lm}\to\xi^{\lm_0}$ in $L^2({\cal F}_T)$ as
$\lm\to\lm_0$, Assumption (H2'), (L) and (H4') hold for $g^\lm$.
Moreover $(y_t^{\lm},z_t^{\lm})_{t\in[0,T]}$ and
$(y_t^{\lm_0},z_t^{\lm_0})_{t\in[0,T]}$ are the solutions of the
BSDEs $(g^\lm,T,\xi^\lm)$ and $(g^{\lm_0},T,\xi^{\lm_0})$
respectively, then
\begin{eqnarray}\label{4.2}
E[\sup_{t\in[0,T]}\abs{y_t^{\lm}-y_t^{\lm_0}}^2]&\le&
CE\abs{\xi^{\lm}-\xi^{\lm_0}}^2\nonumber\\
&+&CE\int_0^T\abs{g^{\lm}(t,y_t^{\lm_0},z_t^{\lm_0})-g^{\lm_0}(t,y_t^{\lm_0},z_t^{\lm_0})}^2\,ds
\end{eqnarray}
where $C$ is a positive constant only depending on Lipschitz
constant $K$. Moreover, we have
\begin{equation}\label{4.3}
\lim_{\lm\to\lm_0}E[\sup_{t\in[0,T]}\abs{y_t^{\lm}-y_t^{\lm_0}}^2]=0.
\end{equation}
\end{lem}
\begin{proof}
With the usual techniques of BSDE we can get inequality
(\ref{4.2})(see [3] for detail). Because of the continuity of
$g^{\lm}$ in $\lm=\lm_0$ and Lebesgue dominated convergence theorem
we take limit to both sides of (\ref{4.2}) and get equation
(\ref{4.3}). The proof is complete.
\end{proof}

Now, we introduce the approximation sequences of $g^\lm$ as
follows:
\begin{equation}\label{4.4}
\ul{g}^{\lm}_m(t,y,z)=\inf_{(u,v)\in
\Real^{1+d}}\set{g^{\lm}(t,u,v)+m(\abs{y-u}+\abs{z-v})},
\end{equation}
and
\begin{equation}\label{4.5}
\bar{g}^{\lm}_m(t,y,z)=\sup_{(u,v)\in
\Real^{1+d}}\set{g^{\lm}(t,u,v)-m(\abs{y-u}+\abs{z-v})}.
\end{equation}
\begin{lem}\label{lem4.2}
If $g^{\lm}$ satisfies (H1')---(H4'), then for any $m\ge A$, we
have

(1). $\abs{\ul{g}^{\lm}_m(t,y,z)}\le A(\abs{y}+\abs{z}+1),$ and
$\abs{\bar{g}^{\lm}_m(t,y,z)}\le A(\abs{y}+\abs{z}+1)$, for any
$y\in\Real,z\in\Real^d,\lm\in D$ and $t\in[0,T]$.

(2). For any given $y\in\Real,z\in\Real^d,\lm\in D$ and
$t\in[0,T]$, $\ul{g}^{\lm}_m(t,y,z)$ is nondecreasing in $m$ and
$\bar{g}^{\lm}_m(t,y,z)$ is non-increasing in $m$.

(3). $\ul{g}^{\lm}_m$ and $\bar{g}^{\lm}_m$ are Lipschitz
continuous in $(y,z)$, that is, for any $y_1,y_2\in\Real$,
$z_1,z_2\in\Real^d$ and $\lm\in D$, we have
$\abs{\ul{g}^{\lm}_m(t,y_1,z_1)-\ul{g}^{\lm}_m(t,y_2,z_2)}\le
m(\abs{y_1-y_2}+\abs{z_1-z_2}),$ and
$\abs{\bar{g}^{\lm}_m(t,y_1,z_1)-\bar{g}^{\lm}_m(t,y_2,z_2)}\le
m(\abs{y_1-y_2}+\abs{z_1-z_2})$.

(4). If $(y_m,z_m)\to(y,z)$ as $m\to\infty$, then
$\ul{g}^{\lm}_m(t,y_m,z_m)\to g^{\lm}(t,y,z)$, and
$\bar{g}^{\lm}_m(t,y_m,z_m)\to g^{\lm}(t,y,z)$ as $m\to\infty$.

(5). Both $\ul{g}^{\lm}_m$ and $\bar{g}^{\lm}_m$ are continuous in
$\lm=\lm_0$.
\end{lem}
\begin{proof}
It is easy to check (1)---(4) (see \cite{LepeltierMartin1997}). Now,
we will prove (5). For any $\eps>0$, by the definition of
$\ul{g}^{\lm}_m$, there exist $(y^{\eps,\lm}, z^{\eps,\lm})$ and
$(y^{\eps,\lm_0}, z^{\eps,\lm_0})$ such that
\begin{eqnarray*}
&&g^{\lm}(t,y^{\eps,\lm},z^{\eps,\lm})+m\abs{y-y^{\eps,\lm}}+\abs{z-z^{\eps,\lm}}-\eps\le\ul{g}^{\lm}_m(t,y,z)\\
&&\le
g^{\lambda}(t,y^{\eps,\lm_0},z^{\eps,\lm_0})+m\abs{y-y^{\eps,\lm_0}}+\abs{z-z^{\eps,\lm_0}}
\end{eqnarray*}
and
\begin{eqnarray*}
&&g^{\lm_0}(t,y^{\eps,\lm_0},z^{\eps,\lm_0})+m\abs{y-y^{\eps,\lm_0}}+\abs{z-z^{\eps,\lm_0}}-\eps\le\ul{g}^{\lm_0}_m(t,y,z)\\
&&\le
g^{\lm_0}(t,y^{\eps,\lm},z^{\eps,\lm})+m\abs{y-y^{\eps,\lm}}+\abs{z-z^{\eps,\lm}}
\end{eqnarray*}
thus
\begin{eqnarray*}
&&g^{\lm}(t,y^{\eps,\lm},z^{\eps,\lm})-g^{\lm_0}(t,y^{\eps,\lm},z^{\eps,\lm})-\eps\\
&&\le\ul{g}^{\lm}_m(t,y,z)-\ul{g}^{\lm_0}_m(t,y,z)\\
&&\le
g^{\lm}(t,y^{\eps,\lm_0},z^{\eps,\lm_0})-g^{\lm_0}(t,y^{\eps,\lm_0},z^{\eps,\lm_0})+\eps
\end{eqnarray*}
Because $g^{\lm}$ is continuous when $\lambda=\lambda_0$ uniformly
with respect to $(y,z)$, we obtain the continuity of
$\ul{g}^{\lm}_m$ and $\bar{g}^{\lm}_m$ in $\lm=\lm_0$. The proof is
complete.
\end{proof}

\begin{lem}\label{lem4.3}
If $g^\lm$ satisfies (H1')---(H4'), and the processes
$(\ul{y}_t^{\lm,m},\ul{z}_t^{\lm,m})_{t\in[0,T]}$ and
$(\bar{y}_t^{\lm,m},\bar{z}_t^{\lm,m})_{t\in[0,T]}$ are the unique
solutions of the BSDEs $(\ul{g}^{\lm}_m,T,\xi^\lm)$ and
$(\bar{g}^{\lm}_m,T,\xi^\lm)$ respectively, then, for any $\lm\in
D$, we have
$$(\ul{y}_t^{\lm,m},\ul{z}_t^{\lm,m})_{t\in[0,T]}\to(\ul{y}_t^{\lm},\ul{z}_t^{\lm})_{t\in[0,T]},\ and\ (\bar{y}_t^{\lm,m},\bar{z}_t^{\lm,m})_{t\in[0,T]}\to(\bar{y}_t^{\lm},\bar{z}_t^{\lm})_{t\in[0,T]},$$
in $\Ss\times\Hh_d$ as $m\to\infty$, where
$(\ul{y}_t^{\lm},\ul{z}_t^{\lm})_{t\in[0,T]}$ and
$(\bar{y}_t^{\lm},\bar{z}_t^{\lm})_{t\in[0,T]}$ are the minimal
solution and maximal solution of BSDE (\ref{4.1}).
\end{lem}
Now, we give our result for the general case.
\begin{theorem}\label{thm4.4}
If $g^\lm$ satisfies (H1')---(H4'), then the following statements
are equivalent:

(iii). Uniqueness: there exists a unique solution of BSDE
(\ref{4.1}) when $\lm=\lm_0$, that is, the solution of
$(g^{\lm_0},T,\xi^{\lm_0})$ is unique.

(iv). Continuous dependence with respect to $g$ and $\xi$: for any
$\xi^{\lm}$, $\xi^{\lm_0}\in L^2({\F}_T)$, if
$\xi^{\lm}\to\xi^{\lm_0}$ in $L^2({\F}_T)$ as $\lm\to\lm_0$,
$(y_t^{\lm},z_t^{\lm})_{t\in [0,T]}$ are any solutions of BSDEs
(\ref{4.1}), $(y_t^{\lm_0},z_t^{\lm_0})_{t\in [0,T]}$ is any
solution of BSDE (\ref{4.1}) when $\lm=\lm_0$, then
$$\lim_{\lm\to\lm_0}E[\sup_{t\in[0,T]}\abs{y_t^{\lm}-y_t^{\lm_0}}^2]=0.$$
\end{theorem}
\begin{proof}
This proof is similar to that of Theorem \ref{thm3.1}. For the
sake of completeness, we give the sketch of proof. Firstly, we
prove (iii) implies (iv). We can get the inequalities similarly to
(\ref{3.6}), that is, $\ul{y}_t^{m,\lm}\le \ul{y}_t^{\lm}\le
y_t^{\lm}\le \bar{y}_t^{\lm}\le \bar{y}_t^{m,\lm}$, for any
$t\in[0,T]$ and $m\ge A$. So,
\begin{eqnarray*}
E[\sup_{t\in[0,T]}\abs{y_t^{\lm}-y_t^{\lm_0}}^2]&\le&
2E[\sup_{t\in[0,T]}\abs{\ul{y}_t^{\lm,m}-\ul{y}_t^{\lm_0,m}}^2]+2E[\sup_{t\in[0,T]}\abs{\ul{y}_t^{\lm_0,m}-y_t^{\lm_0}}^2]\\
&+&2E[\sup_{t\in[0,T]}\abs{\bar{y}_t^{\lm,m}-\bar{y}_t^{\lm_0,m}}^2]+2E[\sup_{t\in[0,T]}\abs{\bar{y}_t^{\lm_0,m}-y_t^{\lm_0}}^2]
\end{eqnarray*}
Fixed $m$, with the help of Lemma \ref{lem 4.1} and Lemma
\ref{lem4.2} and the continuity of $\ul{g}^{\lm}_m$ and
$\bar{g}^{\lm}_m$ when $\lm=\lm_0$, we have,
$$E[\sup_{t\in[0,T]}\abs{\ul{y}_t^{\lm,m}-\ul{y}_t^{\lm_0,m}}^2]\to
0,\ and\
E[\sup_{t\in[0,T]}\abs{\bar{y}_t^{\lm,m}-\bar{y}_t^{\lm_0,m}}^2]\to
0$$
as $\lm\to\lm_0$, for any $m\ge A$. By Lemma \ref{lem4.3} and
the uniqueness of solution for $(g^{\lm_0},T,\xi^{\lm_0})$
(Condition (iii)), we obtain, as $m\to\infty$,
$$E[\sup_{t\in[0,T]}\abs{\ul{y}_t^{\lm_0,m}-y_t^{\lm_0}}^2]\to 0,\
and\ E[\sup_{t\in[0,T]}\abs{\bar{y}_t^{\lm_0,m}-y_t^{\lm_0}}^2]\to
0.$$ This implies (iv).

Now we will prove that (iv) implies (iii). Take
$\xi^{\lm}=\xi^{\lm_0}$, $g^{\lm}=g^{\lm_0}$. For equation
(\ref{4.1}), set $y_t^{\lm}:=\bar{y}_t^{\lm}=\bar{y}_t^{\lm_0}$.
For equation $(g^{\lm_0},T,\xi^{\lm_0})$, take
$y_t^{\lm_0}=\ul{y}_t^{\lm_0}$. By (iv), we have
$\bar{y}_t^{\lm_0}=\ul{y}_t^{\lm_0}$. The proof is complete.
\end{proof}

\end{document}